\documentclass[11pt]{article}
\usepackage{graphics,amssymb,amsmath,epsf}

\usepackage{fancyhdr}

\begin{document}

\title{ Inference on Difference of Means of two Log-Normal Distributions; A Generalized Approach}
\author{K. Abdollahnezhad$^1$,  M. Babanezhad$^2$, A. A. Jafari$^3$\\
{\small $^1$Department of Statistics, Golestan University, Gorgan, Iran}\\
{\small $^2$Department of Statistics, Golestan University, Gorgan, Iran}\\
{\small $^3$Department of Statistics, Yazd University, Yazd, Iran}.
 }
\date{}
\maketitle

\begin{abstract}

Over the past decades, various methods for comparing the means of
two log-normal have been proposed. Some of them are differing in terms
of how the statistic test adjust to accept or to reject the null hypothesis.
In this study, a new method of test for comparing the means of two log-
normal populations is given through the generalized measure of evidence
to have against the null hypothesis. However calculations of this method
are simple, we find analytically that the considered method is doing well
through comparing the size and power statistic test. In addition to the
simulations, an example with real data is illustrated.

\end{abstract}


\section{Introduction}

One often encounters with random variables that are inherently positive in
some real life applications such as analyzing biological, medical, and industrial
data. In this regards the normal distribution is applied in most of applications.
In the family of normal distribution the Log normal distribution has
a long term applications. In probability theory, a log-normal distribution is
a continuous probability distribution of a random variable whose logarithm is
normally distributed. Further, a variable might be modeled as log-normal if it
can be thought of as the multiplicative product of many independent random
variables each of which is positive. The suitability of the log-normal random
variable has been investigated by some researchers (Crow and Shimizu, 1988).
There are also some recent articles regarding the statistical inference of pa-
rameters of several log-normal distributions. For example, one-sided test have
been investigated for two distributions with a large sample under the homo-
geneity of the mean parameters for m log-normal populations (Zhou et al.,
1997; Ahmed et al., 2002). Further, exact confidence interval test for the ratio
or difference of the means of two log-normal distributions using the generalized
variable and generalized p-values through a modified likelihood ratio has been
done (Krishnamoorthy and Mathew, 2003; Gill, 2004; Gupta and Li, 2005). In
this paper, we consider random samples from two lognormal populations and
our interest is to present a test of difference of the means of these two popu-
lations. In Section 2, the theory of generalized p-value is introduced. Section
3 is devoted to an exact one-sided test or two-sided test for two log-normal
distributions. We compare the size and power of different proposed methods
to test of the means of two log-normal populations in Section 4 through simu-
lation. We examine them by a numerical example with real data set. A brief
discussion is given in Section 5.

\section{Generalized $p$-value}

The concept of generalized $p$-value was first introduced by Tsui
and Weerahandi (1989) to deal with the statistical testing
problem in which nuisance parameters are present and it is
difficult or impossible to obtain a nontrivial test with a fixed
level of significance. The setup is as follows. Let
\mbox{\boldmath $X$} be a random variable having density function
$f(\mbox{\boldmath $x$}| \mbox{\boldmath $\zeta$})$, where
$\mbox{\boldmath $\zeta$}=(\theta ,\mbox{\boldmath $\eta$})$ is a
vector of unknown parameters, $\theta $ is the parameter of
interest, and \mbox{\boldmath $\eta$} is a vector of nuisance
parameters. Suppose we are interested to test
\begin{equation}
H_{\circ }:\theta \leqslant \theta _{\circ }\text{ \ \ }vs\text{ \ \ }%
H_{1}:\theta >\theta _{\circ }\text{,}
\end{equation}
where $\theta _{\circ }$ is a specified value.

Let $\mbox{\boldmath $x$}$ denote the observed value of
$\mbox{\boldmath $X$}$ and consider a variable $T(\mbox{\boldmath
$X$};\mbox{\boldmath $x$},\mbox{\boldmath $\zeta$} )$, by the name
of generalized variable. We assume that $T(\mbox{\boldmath
$X$};\mbox{\boldmath $x$},\mbox{\boldmath $\zeta$} )$ satisfies
the following conditions:

(i) For fixed $\mbox{\boldmath $x$}$, the distribution of
$T(\mbox{\boldmath $X$};\mbox{\boldmath $x$},\mbox{\boldmath
$\zeta$} )$ is free from the nuisance parameters $\mbox{\boldmath
$\eta$}$.

(ii) $t_{obs}=T(\mbox{\boldmath $x$};\mbox{\boldmath
$x$},\mbox{\boldmath $\zeta$})$ is free from any unknown
parameters.

(iii) For fixed $\mbox{\boldmath $x$}$ and $\mbox{\boldmath
$\eta$}$, $T(\mbox{\boldmath $X$};\mbox{\boldmath
$x$},\mbox{\boldmath $\zeta$} )$ is either stochastically
increasing or decreasing in $\theta $ for any given $t.$

Under the above conditions, if $T(\mbox{\boldmath
$X$};\mbox{\boldmath $x$},\mbox{\boldmath $\zeta$} )$ is
stochastically increasing in $\theta $, then the generalized
$p$-value for testing the hypothesis in (1) can be defined as
\begin{equation}
p=\sup_{\theta \leqslant \theta _{\circ }}\text{
}P(T(\mbox{\boldmath $X$};\mbox{\boldmath $x$},\theta
,\mbox{\boldmath $\eta$})\geq t^{\ast })=P(T(\mbox{\boldmath
$X$};\mbox{\boldmath $x$},\theta_{\circ} ,\mbox{\boldmath $\eta$})
)\geq t^{\ast }),
\end{equation}
where $t^{\ast }=T(\mbox{\boldmath $X$};\mbox{\boldmath
$x$},\theta_{\circ} ,\mbox{\boldmath $\eta$}) )$.

For further details and for several applications based on the
generalized $p$-value, we refer to the book by Weerahandi (1995).

\section{A generalized test variable}

Let $Y_{ij}=\ln (X_{ij})\thicksim N(\mu _{i},\sigma _{i}^{2})$,
$i=1,2$ , $j=1,2,...,n_{i}$ be independent random samples from
two log-normal populations. We know that $M_{i}=E(X_{ij})=\exp
(\mu _{i}+0.5\sigma _{i}^{2}).$ The problem of our interest  is
one sided and two sided test hypothesis about $\eta =M_{1}-M_{2}
$.

In this section, using the concept of generalized $p$-value, we
test
\begin{equation}
H_{\circ }:M_{1}\leqslant M_{2}\text{ \ }vs\text{ \ \ }H_{1}:M_{1}>M_{2}%
\text{,}
\end{equation}
which is equivalent to
\begin{equation}
H_{\circ }:\theta \text{ }\leqslant 0\text{\quad }vs\text{ \ \ }H_{1}:\theta
\text{ }>0\text{,}
\end{equation}
where $\theta =\ln M_{1}-\ln M_{2}$.

The MLE's for $\mu _{i}$ and $\sigma _{i}^{2}$ $(i=1,2)$ are
$\bar{Y}_{i}$ and $S_{i}^{2}$, respectively, where
\[
\bar{Y}_{i}=\frac{1}{n_{i}}\sum_{i=1}^{n}Y_{ij}\quad ,\quad S_{i}^{2}=\frac{1%
}{n_{i}}\sum_{i=1}^{n}(Y_{ij}-\bar{Y}_{i})^{2}.
\]

Now, consider
\begin{eqnarray*}
T
&=&\bar{y}_{1.}-\bar{y}_{2.}+\dfrac{\bar{Y}_{2.}-\bar{Y}_{1.}-(\mu
_{2}-\mu _{1})}{\sqrt{\dfrac{\sigma
_{1}^{2}}{n_{1}}+\dfrac{\sigma
_{2}^{2}}{n_{2}}}}\sqrt{\dfrac{\sigma
_{1}^{2}s_{1}^{2}}{n_{1}S_{1}^{2}}+\dfrac{\sigma
_{2}^{2}s_{2}^{2}}{n_{2}S_{2}^{2}}}+\dfrac{\sigma
_{1}^{2}s_{1}^{2}}{2S_{1}^{2}}-\dfrac{\sigma _{2}^{2}s_{2}^{2}}{2S_{2}^{2}}-\theta  \\
&=&\bar{y}_{1.}-\bar{y}_{2.}+Z\sqrt{\dfrac{s_{1}^{2}}{U_{1}}+\dfrac{s_{2}^{2}
}{U_{2}}}+\dfrac{n_{1}s_{1}^{2}}{2U_{1}}-\dfrac{n_{2}s_{2}^{2}}{2U_{2}}
-\theta ,
\end{eqnarray*}
where
\[
Z=\dfrac{\bar{Y}_{2.}-\bar{Y}_{1.}-(\mu _{2}-\mu
_{1})}{\sqrt{\dfrac{\sigma _{1}^{2}}{n_{1}}+\dfrac{\sigma
_{2}^{2}}{n_{2}}}}\thicksim N(0,1),
\]
and
\[
U_{i}=\frac{n_{i}S_{i}^{2}}{\sigma _{i}^{2}}\thicksim \chi _{(n_{i}-1)}^{2}%
\text{ , \ }i=1,2,
\]
are three independent random variables, and $\bar{y}_{i}$ and
$s_{i}^{2}$ are observed values of $\bar{Y}_{i}$ and $S_{i}^{2}$,
respectively$.$Then, $T$ is a generalized variable for $\theta $
because

i) $\ t_{obs}=0$

ii) distribution of $T$ is free from the nuisance parameters $\mu _{i}$ and $%
\sigma _{i}^{2}$.

iii) the distribution of $T$ is an increasing function with respect to $%
\theta $.

\bigskip

Thus the generalized $p$-value for the null hypothesis (3) is
given by
\begin{equation}
p=P(T\leqslant t_{obs}|\theta =0)=E(\Phi (\dfrac{\bar{y}_{2.}-\bar{y}%
_{1.}+\dfrac{n_{2}s_{2}^{2}}{2U_{2}}-\dfrac{n_{1}s_{1}^{2}}{2U_{1}}}{%
\sqrt{\dfrac{s_{1}^{2}}{U_{1}}+\dfrac{s_{2}^{2}}{U_{2}}}})),
\end{equation}
where $\Phi (.)$ is the standard normal distribution function and the
expectation is taken with respect to independent chi-square random
variables, $U_{1}$ and $U_{2}.$

This generalized $p$-value can be well approximated by a Monte
Carlo simulation using the following algorithm:

\bigskip

{\bf Algorithm 1.} For a given data set $x_{i1},...,x_{in_{i}},$ set $%
y_{ij}=\ln (x_{ij})$, $i=1,...,k$ , $j=1,2.$

Compute $\bar{y}_{1.},$ $\bar{y}_{2.},$ $s_{1}^{2},$ $s_{2}^{2}$

For $l=1$ to $m$

Generate $U_{1}\thicksim \chi _{(n_{1}-1)}^{2}$ , $U_{2}\thicksim \chi
_{(n_{2}-1)}^{2}.$

Calculate $T_{l}=\Phi (\dfrac{\bar{y}_{2.}-\bar{y}_{1.}+\dfrac{
n_{2}s_{2}^{2}}{2U_{2}}-\dfrac{n_{1}s_{1}^{2}}{2U_{1}}}{\sqrt{\dfrac{%
s_{1}^{2}}{U_{1}}+\dfrac{s_{2}^{2}}{U_{2}}}}).$

\noindent $\dfrac{1}{m}\sum\limits_{l=1}^{m}T_{l}$ is a Monte
Carlo estimate of generalized $p$-value for the null hypothesis
(3).

\bigskip

The generalized $p$-value in (5) is used for one sided test
hypothesis but we can use this generalized $p$-value for two sided
test hypothesis by
\[
p=2\min \{p,1-p\},
\]
where $p$ is the generalized $p$-value in (5).

\section{ Simulation Study}

To investigate the power of the considered test statistics in finite samples, we
conducted a simulation experiment. To do so, several data set from two log-
normal distributions with $\mu_2 = 0$ were generated. For each scenarios 10000
sample size are performed. The size and the power of the considered test statis-
tics are summarized in tables in table \ref{table2} and \ref{table3}. These tests are (a) generalized
p-value in (5) (b) generalized p-value by Krishnamoorthy and Mathew (2003)
(c) Z-score test by Zhou et al. (1997). The simulation study indicates that
(i) The size for (a) and (b) are close to 0.05 and the powers are close to each
other.
(ii) The size of (c) is very larger than nominal level, 0.05.

\section{ Numerical examples}

The data show the amount of rainfall (in acre-feet) from 52
clouds; 26 clouds were chosen at random and seeded with silver
nitrate. We can show that log-normal model fits the data. The
summary statistics for the log-transformed data are given in Table
\ref{table1}.

\begin{table}[h]
\begin{center}
\caption{ The summary statistics for the log-transformed data of
rainfall } \label{table1}

\begin{tabular}{lccc}
\hline\hline Clouds & $n_{i}$ & $\bar{y}_{i.}$ & $s_{i}^{2}$ \\
\hline
\begin{tabular}{l}
seeded clouds \\
unseeded clouds
\end{tabular}
&
\begin{tabular}{l}
26 \\
26
\end{tabular}
&
\begin{tabular}{l}
5.134 \\
3.990
\end{tabular}
&
\begin{tabular}{l}
2.46 \\
2.60
\end{tabular}
\\ \hline
\end{tabular}

\end{center}
\end{table}

In order to understand the effect of silver nitrate seeding, we
like to test

\begin{equation}
H_{\circ }:M_{1}=M _{2}\text{ \ }vs\text{ \ \ }H_{1}:M _{1}>M
_{2},
\end{equation}
where $M _{i}=\exp (\mu _{i}+0.5\sigma _{i}^{2}),$ $i=1,2.$

The $p$-values for our generalized approach,  Krishnamoorthy and
Mathew approah  and Z-score test  are 0.0779, 0.0747 and 0.0599
respectively. Therefore, we cannot reject $H_{\circ }$ at the
level of $0.05$, using all 3 methods.

\begin{table}[h]
\begin{center}
\caption{Simulated sizes of the tests at 5\% significance level
when $\mu _{2}=0$.} \label{table2}

\begin{tabular}{cccc}
\hline\hline
\begin{tabular}{r}
$n_{1}$%
\end{tabular}
\
\begin{tabular}{r}
$n_{2}$%
\end{tabular}
\ \ \ \
\begin{tabular}{r}
$\mu _{1}$%
\end{tabular}
\begin{tabular}{r}
$\sigma _{1}^{2}$%
\end{tabular}
\
\begin{tabular}{r}
$\sigma _{2}^{2}$%
\end{tabular}
& (a) & (b) & (c) \\ \hline\hline
\begin{tabular}{r}
4 \\
\\
\\
\\
10 \\
\\
\\
\\
25 \\
\\
\\
\\
\\
\\
40 \\
\\
\\
25 \\
\\
\\
40 \\
25 \\
40 \\
\\
100 \\
\\
\\
25
\end{tabular}
\begin{tabular}{r}
4 \\
\\
\\
\\
10 \\
\\
\\
\\
25 \\
\\
\\
\\
\\
\\
25 \\
\\
\\
40 \\
\\
\\
25 \\
40 \\
40 \\
\\
25 \\
\\
\\
100
\end{tabular}
\ \
\begin{tabular}{r}
1 \  \\
0 \\
5 \\
0 \\
1 \\
0 \\
5 \\
0 \\
0 \\
0 \\
0 \\
0 \\
2 \\
4 \\
0 \\
0 \\
0 \\
0 \\
0 \\
0 \\
5 \\
5 \\
8 \\
14 \\
0 \\
0 \\
0 \\
0
\end{tabular}
\begin{tabular}{r}
2 \\
3 \\
2 \\
12 \\
2 \\
3 \\
2 \\
12 \\
1 \\
5 \\
10 \\
100 \\
4 \\
8 \\
1 \\
5 \\
10 \\
1 \\
5 \\
10 \\
2 \\
2 \\
4 \\
4 \\
1 \\
5 \\
10 \\
1
\end{tabular}
\begin{tabular}{r}
4 \\
3 \\
12 \\
12 \\
4 \\
3 \\
12 \\
12 \\
1 \\
5 \\
10 \\
100 \\
8 \\
16 \\
1 \\
5 \\
10 \\
1 \\
5 \\
10 \\
12 \\
12 \\
20 \\
32 \\
1 \\
5 \\
10 \\
1
\end{tabular}
& \ \ \
\begin{tabular}{r}
421 \\
344 \\
464 \\
392 \\
612 \\
546 \\
515 \\
538 \\
512 \\
521 \\
486 \\
521 \\
538 \\
492 \\
391 \\
412 \\
459 \\
382 \\
394 \\
435 \\
521 \\
312 \\
536 \\
513 \\
451 \\
374 \\
396 \\
482
\end{tabular}
&
\begin{tabular}{r}
436 \\
405 \\
510 \\
391 \\
603 \\
581 \\
552 \\
538 \\
524 \\
522 \\
531 \\
531 \\
520 \\
493 \\
467 \\
425 \\
416 \\
376 \\
373 \\
412 \\
510 \\
341 \\
492 \\
546 \\
473 \\
379 \\
388 \\
464
\end{tabular}
&
\begin{tabular}{r}
1091 \\
367 \\
2168 \\
112 \\
895 \\
432 \\
1433 \\
386 \\
614 \\
516 \\
446 \\
396 \\
828 \\
851 \\
506 \\
491 \\
512 \\
364 \\
244 \\
199 \\
1061 \\
586 \\
932 \\
922 \\
664 \\
714 \\
720 \\
295
\end{tabular}
\\ \hline
\end{tabular}

\end{center}

\end{table}

\begin{table}[]
\begin{center}
\caption{Simulated powers of the tests at 5\% significance level
when $\mu _{2}=0$.} \label{table3}

\begin{tabular}{cccc}
\hline\hline
\begin{tabular}{r}
$\ \ \ \ n_{1}\ \ \ \ n_{2\text{\ \ \ \ }}$ $\ \mu _{1}$\ $\ \sigma _{1}^{2}$%
\ $\ \ \sigma _{2}^{2}$%
\end{tabular}
& (a) & (b) & (c) \\ \hline\hline
\begin{tabular}{r}
4 \\
\\
\\
\\
10 \\
\\
\\
\\
25 \\
\\
\\
\\
\\
\\
40 \\
\\
\\
25 \\
\\
\\
40 \\
\\
25 \\
\\
100 \\
\\
\\
25
\end{tabular}
\begin{tabular}{r}
4 \\
\\
\\
\\
10 \\
\\
\\
\\
25 \\
\\
\\
\\
\\
\\
25 \\
\\
\\
40 \\
\\
\\
25 \\
\\
40 \\
\\
25 \\
\\
\\
100
\end{tabular}
\ \
\begin{tabular}{r}
0\  \\
3 \\
0 \\
4 \\
0 \\
0 \\
3 \\
4 \\
1 \\
1 \\
0 \\
1 \\
0 \\
0 \\
1 \\
1 \\
1 \\
1 \\
1 \\
1 \\
1 \\
1 \\
1 \\
1 \\
1 \\
1 \\
1 \\
1
\end{tabular}
\begin{tabular}{r}
12 \\
2 \\
20 \\
1 \\
12 \\
20 \\
2 \\
1 \\
1 \\
5 \\
4 \\
10 \\
9 \\
4 \\
1 \\
5 \\
10 \\
1 \\
5 \\
10 \\
5 \\
10 \\
5 \\
10 \\
1 \\
5 \\
10 \\
1
\end{tabular}
\begin{tabular}{r}
4 \\
4 \\
4 \\
1 \\
4 \\
4 \\
4 \\
1 \\
1 \\
5 \\
2 \\
10 \\
7 \\
1 \\
1 \\
5 \\
10 \\
1 \\
5 \\
10 \\
4 \\
9 \\
4 \\
9 \\
1 \\
5 \\
10 \\
1
\end{tabular}
& \ \ \
\begin{tabular}{r}
1523 \\
1261 \\
2610 \\
5753 \\
4136 \\
6941 \\
2961 \\
9931 \\
8370 \\
1916 \\
3564 \\
1126 \\
1314 \\
7411 \\
8392 \\
2023 \\
1173 \\
9194 \\
2243 \\
1120 \\
8836 \\
1831 \\
4464 \\
1956 \\
8834 \\
1856 \\
942 \\
9984
\end{tabular}
&
\begin{tabular}{r}
1496 \\
1204 \\
2601 \\
5772 \\
4089 \\
6903 \\
3114 \\
9942 \\
8345 \\
1843 \\
3521 \\
1123 \\
1360 \\
7390 \\
8324 \\
2025 \\
1123 \\
9135 \\
2307 \\
1145 \\
8814 \\
1734 \\
4482 \\
1893 \\
8762 \\
1932 \\
913 \\
9913
\end{tabular}
&
\begin{tabular}{r}
364 \\
3832 \\
334 \\
9621 \\
2334 \\
4562 \\
5173 \\
9990 \\
8917 \\
2081 \\
3157 \\
1250 \\
1225 \\
6854 \\
9036 \\
2736 \\
1492 \\
9401 \\
2159 \\
784 \\
4649 \\
2263 \\
4955 \\
1493 \\
9512 \\
3364 \\
1825 \\
9893
\end{tabular}
\\ \hline
\end{tabular}

\end{center}
\end{table}


\begin{thebibliography}{}

\bibitem{}Ahmed, S. E., Tomkins, R. J. and Volodin, A.I. (2001). Test of
homogeneity of parallel samples from lognormal populations with
unequal variances, {\it Journal of Statistical Research}, 35, no
2, 25-33.

\bibitem{} Crow, E. L. and Shimizu, K. (1988). {\it Lognormal distribution},
Marcel Dekker: New York.

\bibitem{} Gill, P. S. (2004). Small sample inference for the comparison
of means of lognormal distribution, {\it Biometrics}, 60, 525-527.

\bibitem{} Gupta, R. C. and Li, X. (2005). Statistical inferences on the common
mean of two log-normal distributions and some applications in
reliability, appeared in {\it Computational Statistics and Data
Analysis}.

\bibitem{} Krishnamoorthy, K. and Mathew, T. (2003). Inferences on the means of
lognormal distributions using generalized p-values and generalized
confidence interval, {\it Journal of Statistical Planning and
Inference}, 115, 103-121.

\bibitem{} Krishnamoorthy, K. and Yong Lu. (2003). Inferences on the common mean
of several normal populations based on the generalized variable method, {\it %
Biometrics}, 59, 237-247.

\bibitem{} Tsui, K. W. and Weerahandi, S. (1989). Generalized p-values in
significance testing of hypothesis in the presence of nuisance
parameters, {\it J. Am. Statist. Assoc}., 84, 602-607.

\bibitem{}  Weerahandi, S. (1993). Generalized confidence intervals, {\it J. Am.
Statist. Assoc}., 88, 899-905.

\bibitem{} Weerahandi, S. (1995a). {\it Exact statistical methods for data
analysis}, Springer, NewYork.

\bibitem{} Weerahandi, S. and Berger, V. W. (1999). Exact inference for growth
curves with interclass correlation structure, {\it Biometrics},
55, 921-924.

\bibitem{} Zhou, X. H., Gao, S. and Hui. S.L. (1997). Methods for comparing the
means of two independent lognormal samples, {\it Biometrics}, 53,
1129-1135.

\bibitem{} Zhou, X. H. and Tu. W. (1999). Comparison of several
independent population means when their samples contain lognormal
and possibly zero observations, {\it Biometrics}, 55, 645-651.


\end{thebibliography}
\end{document}